\documentclass{amsart}
\usepackage{palatino}
\usepackage{amsthm, amsmath}
\usepackage{amssymb} 
\usepackage{mathrsfs}
\usepackage[hmargin=3cm, vmargin=2.4cm]{geometry}
\usepackage[breaklinks,bookmarksopen,bookmarksnumbered]{hyperref}

\theoremstyle{plain}
\newtheorem*{thm}{Theorem}
\newtheorem{lem}{Lemma}

\newtheorem{cor}{Corollary}

\newcommand*{\lr}[1]{\langle {#1}\rangle}
\newcommand\pr{\noindent\textit{Proof} : }

\newcommand\Ker{\operatorname{Ker}}
\newcommand{\mo}{\smallsetminus}

\newcommand\Z{\mathbb{Z}}
\newcommand\Q{\mathbb{Q}}

\renewcommand\P{\mathbb{P}}
\renewcommand\O{\mathcal{O}}

\begin{document}
\title{A non-hyperelliptic curve with  torsion Ceresa cycle modulo algebraic equivalence}
\author{Arnaud Beauville}
\address{Universit\'e C\^ote d'Azur\\
CNRS -- Laboratoire J.-A. Dieudonn\'e\\
Parc Valrose\\
F-06108 Nice cedex 2, France}
\email{arnaud.beauville@unice.fr}
 \author{Chad Schoen}
 \address{Department of Mathematics, Duke University Box 90320, Durham NC 27708-0320, USA}
 \email{schoen@math.duke.edu}
 \thanks{Thanks to E. Colombo and B. van Geemen for their crucial input.
The second named author thanks S. Katz and M. Reid for helpful discussions.}
\begin{abstract}
We exhibit a non-hyperelliptic  curve $C$ of genus $3$ such that the class of the Ceresa cycle $[C]-[C^-]$ in $JC$ modulo algebraic equivalence is torsion.
\end{abstract}
\maketitle 
\section{Introduction}
Let $C$ be a complex curve  of genus $g\geq 2$, and $p$ a point of $C$. We embed $C$ into its Jacobian $J$ by the Abel-Jacobi map $x\mapsto [x]-[p]$; we denote by $C^{-}$ the image of $C$ under the involution $(-1)_J: a\mapsto -a$ of $J$. The 
 \emph{Ceresa class}  is the class $\mathfrak{z}(C):=[C] - [C^{-}]$ in the group   
 $A_1(J)$ of $1$-cycles on $J$ modulo algebraic equivalence (it is independent of the choice of $p$).  Since $(-1)_J$ acts trivially on $H^{p}(J,\Z)$ for $p$ even,  $\mathfrak{z}(C)$ belongs to the \emph{Griffiths group} $G(J)$, the kernel of the cycle class map $A_1(J)\rightarrow H^{2g-2}(J,\Z)$.
 
 Ceresa classes have played a prominent role in the study of Griffiths groups,
especially in the development of techniques for showing that a given element
is non-zero \cite{C,C-P,H}. In addition they played
an important role in showing that $G(J)$ can have infinite rank \cite{N}.
As the conjectures of Bloch and Beilinson were developed and are studied
$\mathfrak{z}(C)$ appears repeatedly \cite{Bl, BST}, \cite[\S1.5]{Z},
always as an element of infinite order.

When $C$ is hyperelliptic, $\mathfrak{z}(C)=0$; in fact $C-C^{-}$ is zero
as a cycle when $p$ is a Weierstrass point.
  In this note we will exhibit what we believe to be the first example of a non-hyperelliptic  curve $C$  with $\mathfrak{z}(C)=0$ in $A_1(J)\otimes \Q$. The curve $C$ has genus $3$, and admits an automorphism $\sigma $ of order 9, such that the quotient variety $J/\lr{\sigma }$ is uniruled. This implies that  the Griffiths group of a resolution of $J/\lr{\sigma }$ is torsion; going back to $J$ gives the result.
 
 \medskip	
 \section{Main result}
\begin{thm}
Let $C\subset \P^2$ be the genus $3$ curve defined by $X^4+XZ^3+Y^3Z=0$. Then $\mathfrak{z}(C)=0$ in $A_1(J)\otimes \Q$.
\end{thm}
\pr Let $\zeta   $ be a primitive 9-th root of unity. We consider the automorphism $\sigma $ of $C$ defined by $\sigma (X,Y,Z)=(X,\zeta ^2Y,\zeta ^3Z)$. We use the fixed point $p=(0,0,1)$ to embed $C$ in its Jacobian $J$, so that the action of $\sigma $ on $J$ preserves $C$ and $C^{-}$. We denote by $V$ the quotient variety $J/\lr{\sigma }$, and by $\pi : J\rightarrow V$ the quotient map. Let $F\subset J$ be the subset of elements with nontrivial stabilizer; the singular locus $\operatorname{Sing}V $ of $V$ is $\pi (F)$. We put $J^{\mathrm{o}}:=J\mo F$ and $V^{\mathrm{o}}:=V\mo \operatorname{Sing}V $.

 \begin{lem}\label{can}
$\operatorname{Sing}V $ is finite; the points $\pi (x)$ for $x\in\Ker(1_J-\sigma )$ 
are non-canonical singularities. 
\end{lem}
\pr 
The space $T_0(J)$ is canonically identified with $H^0(C,K_C)^*$.  
 The elements of $H^0(C,K_C)$ are 
 of  the form  $ L\,\dfrac{XdZ-ZdX}{Y^2Z} $, with $L\in H^0(\P^2,\O_{\P}(1))$ \cite[\S 9.3, Corollary of Theorem 1]{B-K}. It follows that the eigenvalues of $\sigma $ on $H^0(C,K_C)$
 are $ \zeta^5 ,\zeta ^7, \zeta ^8$, and those on $T_0(J)=H^0(C,K_C)^*$ are $\zeta ,\zeta ^2,\zeta ^4$. Therefore $\Ker(1_J-\sigma^{d} )$ is finite for $0<d<9$, so $F$ is finite.
 Since $1+2+4<9$, Reid's criterion \cite[Theorem 3.1]{R} implies that the  singular points $\pi (x)$ for $x\in \Ker(1_J-\sigma )$ are not canonical.\qed
 
 \medskip	
 \begin{lem}
The variety $V$ is uniruled.
\end{lem}
 \pr Let $\rho :\tilde{V}\rightarrow V $ be a resolution of singularities; it suffices to prove that  $\tilde{V} $  has Kodaira dimension $-\infty$ \cite{Mi}. Suppose this is not the case: there exist an integer $r\geq 1$ and a nonzero section $\tilde{\omega}  $ of $K_{\tilde{V}}^r$.   By restriction to $\rho ^{-1}(V^{\mathrm{o}})\cong V^{\mathrm{o}}$, we get  a section $\omega$ of $K^r_{V^{\mathrm{o}}}$, whose pull back under $\pi $ is a nonzero section of $K^r_{J^{\mathrm{o}}}$; therefore $\omega$ is a generator of $K^r_{V^{\mathrm{o}}}$, hence extends to a generator of the reflexive sheaf $K_{V}^{[r]}$ (with the notation of \cite{R}). By construction this generator remains regular on $\tilde{V} $, which means that the singularities of $V$ are canonical \cite[Proposition 1.2]{R}, contradicting Lemma \ref{can}.\qed
 
   \begin{lem}\label{uni}
Let $X$ be a uniruled smooth projective threefold. The Griffiths group $G(X)$ is torsion.
\end{lem}
\pr  
There exists a smooth projective surface $S$ and a dominant rational map $S\times \P^1\dasharrow X$. After blowing up some points and some smooth curves in $S\times \P^1$, we get a smooth projective threefold $W$ birational to $S\times \P^1$, and a
generically finite morphism $f:W\rightarrow X$. Since the Griffiths group is a stably birational invariant (see \cite[Proposition 2.30]{V}), we have $G(W)=G(S)=0$. For $z\in G(X)$, we have $(\deg f)z= f_*f^*(z)=0$, hence $G(X)$ is annihilated by $\deg f$.\qed

\smallskip	
\noindent\emph{Remark}$.-$  One can actually deduce from \cite[Theorem 1 (ii)]{B-S} that $G(X)=0$ -- but  we will not need this fact.

 \medskip	
\noindent\emph{Proof of the Theorem}: We can choose the resolution $\rho :\tilde{V}\rightarrow V $ so that $E:=\rho ^{-1}(\operatorname{Sing V} )$ is a normal crossing divisor, whose irreducible  components are smooth and \emph{rational} \cite[Corollary of Theorem 1]{Fuj}. 
 
 Let $\bar{C}$ and $\bar{C}^-$ be the images in $V$ of $C$ and $C^-$, and 
   let $\tilde{C}  $ and $\tilde{C} ^-$ be their proper transforms in $\tilde{V} $. We have $[\bar{C}]-[\bar{C}^-]=\frac{1}{9} \pi _* ([C]-[C^-])=0$  in $H^4(V^{\mathrm{o}},\Q)$. Now we have an exact sequence \cite[Corollaire 8.2.8]{D}
\[H^2(\tilde{E},\Q )\xrightarrow{\ i_*\ }H^4(\tilde{V},\Q )\rightarrow H^4(V^{\mathrm{o}},\Q)\,,\]where $\tilde{E} $ is the normalization of $E$, and $i$ the composition $\tilde{E}\rightarrow E\hookrightarrow \tilde{V}  $.
Therefore we have $[\tilde{C} ]-[\tilde{C}^- ]=i_*z$ in $H^4(\tilde{V},\Q )$ for some class $z\in H^2(\tilde{E},\Q )$. Since the components of $\tilde{E} $ are rational, $z$ is the class of an element $\mathfrak{z}$ of $A_1(\tilde{E} )\otimes \Q$.
Then $[\tilde{C} ]-[\tilde{C}^- ]-i_*\mathfrak{z}\in A_1(\tilde{V} )\otimes \Q$ is homologous to zero, hence equal to zero by Lemma \ref{uni}. Restricting to $\tilde{V}\mo E \cong V^{\mathrm{o}}$, we get 
 $[\bar{C}]-[\bar{C}^- ]=0$  in $A_1(V^{\mathrm{o}}\otimes \Q)$, hence $[C]-[C^-]=\pi ^*([\bar{C}]-[\bar{C}^- ])=0$  in $A_1(J^{\mathrm{o}})\otimes \Q$. But the restriction map  $A_1(J)\rightarrow A_1(J^{\mathrm{o}})$ is an isomorphism \cite[Example 10.3.4]{Ful}, hence the Theorem.\qed

\medskip	
\section{Complements}

\begin{cor}
Let $\Theta $ be a Theta divisor on $J$. We have $[C]=\dfrac{[\Theta] ^2}{2} $ in $A_1(J)\otimes \Q$ (Poincar\'e formula). 
\end{cor}
\pr Indeed for \emph{any} genus 3 curve $C$ we have $[\Theta] ^2=[C]+[C^-]$ in $A_1(J)\otimes \Q$ (if $p,q$ are two distinct points of $C$, the intersection of $\Theta $ with its translate by $[p]-[q]$ is the union of a translate of $C$ and a translate of $C^{-}$ --- see for instance  \cite[Lecture IV]{Mu}). Thus the corollary is equivalent to the theorem.\qed

\medskip	
Recall that the \emph{modified diagonal cycle} $\Gamma (C,p)$, first considered in \cite{GS}, is the element $\Gamma (C,p)$ of $A_1(C^3)$ defined as follows. We  denote by $[x, x, x], [x, x, p], [x, p, p]$ etc. the classes in $A_1(C\times C\times C)$ of the image of $C$ by the maps $x\mapsto (x,x,x)$, $x\mapsto (x,x,p)$, $x\mapsto (x,p,p)$ etc. Then:
\[\Gamma (C,p):=[x, x, x]-[x, x, p]-[x, p, x]-[p, x, x]+[x, p, p]+[p, x, p]+[p, p, x]\,.\]
By \cite[Remark 3.4]{FLV}, we have
\begin{cor}
$\Gamma (C,p)=0$ in $A_1(C^3)\otimes \Q$.
\end{cor}
\medskip	
Finally let us mention the result of \cite{B3}:  the class of $[C]-[C^-]$ in the intermediate Jacobian $\mathscr{J}_1(J)$ is torsion. It can be also deduced from our theorem, though the proof in \cite{B3} is more direct. 

In \cite{Litt} the authors construct a genus 7 curve with the same property, and suggest that the corresponding Ceresa cycle should be torsion modulo algebraic equivalence (Remark 1.2).


\bigskip	
\begin{thebibliography}{C-vG}

\bibitem[B1]{B1} A. Beauville\,: \textsl{Sur l'anneau de Chow d'une vari\'et\'e ab\'elienne}. Math. Ann. \textbf{273}  (1986), 647-651. 

\bibitem[B2]{B2} A. Beauville\,: \textsl{Algebraic cycles on Jacobian varieties}. Compositio Math. \textbf{140}  (2004), 683-688.

\bibitem[B3]{B3} A. Beauville\,: \textsl{A non-hyperelliptic curve with torsion Ceresa class}. C.R. Acad. Sci. Paris, to appear; preprint \texttt{arXiv:2105.07160}.

\bibitem[BLLS]{Litt}  D. Bisogno, W. Li, D. Litt, P. Srinivasan\,:  \textsl{Group-theoretic Johnson classes and non-hyperelliptic curves with torsion Ceresa class}.  Preprint \texttt{arXiv:2004.06146}. 

\bibitem[Bl]{Bl}   S. Bloch\,: \textsl{Algebraic cycles and values of L-functions}. J. reine angew. Math. \textbf{350}  (1984), 94-107.

\bibitem[Bl-S]{B-S}  S. Bloch, V. Srinivas\,:  \textsl{Remarks on correspondences and algebraic cycles}. Amer. J. Math. \textbf{105}  (1983), no. 5, 1235-1253.

\bibitem[B-K]{B-K}  E. Brieskorn, H. Kn\"orrer\,: \textsl{Plane algebraic curves}. Birkh\"auser Verlag, Basel, 1986.

\bibitem[BST]{BST} J. Buhler, C. Schoen, J. Top\,: \textsl{Cycles, L-functions and triple products of elliptic curves}.
J. Reine Angew. Math. \textbf{492}  (1997), 93-133.


\bibitem[C]{C} G. Ceresa\,:  \textsl{$C$ is not algebraically equivalent to $C^{-}$ in its Jacobian}. Ann. of Math. \textbf{117} (1983), no. 2, 285-291.

\bibitem[C-vG]{CG}  E. Colombo, B. van Geemen\,:  \textsl{Note on curves in a Jacobian}. Compositio Math. \textbf{88}  (1993), 333-353.

\bibitem[C-P]{C-P} A. Collino, G. P. Pirola\,:  \textsl{The Griffiths infinitesimal invariant for a curve in its Jacobian}. Duke Math J. \textbf{78}  (1995), 59-88.

 
\bibitem[D]{D}  P. Deligne\,:  \textsl{Th\'eorie de Hodge, III}. Publ. Math. IHES \textbf{44}  (1974), 5-77.
 
 \bibitem[Fuj]{Fuj}  A. Fujiki\,:  \textsl{On resolutions of cyclic quotient singularities}. Publ. Res. Inst. Math. Sci. \textbf{10}  (1974/75), no. 1, 293-328. 

\bibitem[Ful]{Ful}  W. Fulton\,:  \textsl{Intersection theory}. Ergebnisse der Math. \textbf{2}. Springer-Verlag, Berlin, 1984. 

\bibitem[FLV]{FLV} L. Fu, R. Laterveer, C. Vial\,:  \textsl{Multiplicative Chow-K\"unneth decompositions and varieties of cohomological K3 type}. 
    Ann. Mat. Pura Appl. \textbf{200}  (2021), no. 5, p. 2085-2126. 


\bibitem[G-S]{GS}  B. Gross, C. Schoen\,:  \textsl{The modified diagonal cycle on the triple product of a pointed curve}. Ann. Inst. Fourier \textbf{45}  (1995), 649-679.

\bibitem[H]{H} B. Harris\,:  \textsl{Homological versus algebraic equivalence in a Jacobian}. Proc. Nat. Acad. Sci. U.S.A. \textbf{80}  (1983), no. 4,  1157-1158.


\bibitem[M]{Mi} Y.       Miyaoka\,:  \textsl{On the Kodaira dimension of minimal threefolds}. Math. Ann. \textbf{281} (1988), 325-332.

\bibitem[Mu]{Mu} D. Mumford\,: \textsl{Curves and their Jacobians}.  The red book of varieties and schemes,
Lecture Notes in Math. \textbf{1358}. Springer-Verlag, Berlin, 1999.

\bibitem[N]{N} M. Nori\,: \textsl{Cycles on the generic abelian threefold}.
Proc. Indian Acad. Sci. Math. Sci. \textbf{99}  (1989), no. 3, 191-196.

\bibitem[R]{R} M. Reid\,: \textsl{Canonical $3$-folds}. Journ\'ees de G\'eometrie Alg\'ebrique d'Angers,  pp. 273-310, Sijthoff \& Noordhoff, 1980.
 
\bibitem[V]{V}   C. Voisin\,:  \textsl{Birational invariants and decomposition of the diagonal}. Birational Geometry of Hypersurfaces, pp. 3-71; Lecture Notes of the UMI \textbf{26}, Springer, 2019.  

\bibitem[Z]{Z} S-W Zhang\,: \textsl{Gross-Schoen cycles and dualising sheaves}. Invent. math. \textbf{179} (2010), no. 1, 1-73.

\end{thebibliography}
\end{document}